\documentclass[11pt,a4paper]{amsart}
\usepackage{amsmath,amsthm,amsfonts,amssymb,bbm}
\usepackage{graphicx,psfrag,subfigure,color}
\usepackage{cite,stackrel}
\usepackage{hyperref}

\numberwithin{equation}{section}

\newtheorem{prop}{Proposition}[section]
\newtheorem{thm}[prop]{Theorem}

\newtheorem{conj}[prop]{Conjecture}
\newtheorem{cla}[prop]{Claim}
\newtheorem{op}[prop]{Open problem}

\newtheorem{rem}[prop]{Remark}

\title[$(2+1)$-dimensional interface dynamics]{ $(2+1)$-dimensional interface dynamics:  mixing time, hydrodynamic limit and
  Anisotropic KPZ growth} \author{ Fabio Toninelli} \address{CNRS and
  Institut Camille Jordan, Universit\'e Lyon 1, 43 bd du 11 novembre
  1918, 69622 Villeurbanne, France. E-mail: {\tt
    toninelli@math.univ-lyon1.fr}} \date{}

\begin{document}

\maketitle

\begin{abstract}
  Stochastic interface dynamics serve as mathematical models for diverse
  time-dependent physical phenomena: the evolution of boundaries
  between thermodynamic phases, crystal growth, random deposition...
  Interesting limits arise at large space-time scales: after suitable
  rescaling, the randomly evolving interface converges to the solution
  of a deterministic PDE (hydrodynamic limit) and the fluctuation
  process to a (in general non-Gaussian) limit process.  In contrast
  with the case of $(1+1)$-dimensional models, there are very few mathematical
  results in dimension $(d+1), d\ge2$. As far as growth models are
  concerned, the $(2+1)$-dimensional case is particularly interesting:
  Wolf \cite{Wolf} conjectured the existence of two different
  universality classes (called KPZ and Anisotropic KPZ), with
  different scaling exponents.  Here, we review recent mathematical
  results on (both reversible and irreversible) dynamics of some
  $(2+1)$-dimensional discrete interfaces, mostly defined through a
  mapping to two-dimensional dimer models.  In particular, in the
  irreversible case, we discuss mathematical support and remaining
  open problems concerning Wolf's conjecture \cite{Wolf} on the
  relation between the Hessian of the growth velocity on one side, and
  the universality class of the model on the other.
\end{abstract}

\section{Introduction}

Many phenomena in nature involve the evolution of interfaces. A first
example is related to  phenomena of deposition on a substrate, in which case the
interface is the boundary of the deposed material: think for instance
of crystal growth by molecular beam epitaxy or, closer to everyday
experience, of the growth of a layer of snow during snowfall (see
e.g.  \cite{BarSta} for a physicist's introduction to 
growth phenomena). Another example is the evolution of the boundary
between two thermodynamic phases of matter. Think of a block of ice
immersed in water: the shape of the ice block, hence the water/ice
boundary,  changes with time and of course  the dynamics  is very
different according to whether temperature is above, below or exactly
at $0$ $^\circ$C.

A common feature of these examples is that on macroscopic (i.e. large)
scales the interface evolution appears to be deterministic, while a
closer look reveals that the interface is actually rough and presents
seemingly random fluctuations (this is particularly evident in the
snow example, since snowflakes have a visible size).

To try to model mathematically such phenomena, a series of
simplifications are adopted.  First, the so-called \emph{effective
  interface approximation}: the $d$-dimensional interface in
$(d+1)$-dimensional space is modeled as a height function
$h: x\in \mathbb R^d\mapsto h_x(t)\in \mathbb R$, where $h_x(t)$ gives
the height of the interface above point $x$ at time $t$ (think of
$d=2$ in the case of snow falling on your garden, but $d=1$ for
instance for snow falling and sliding down on your car window).  This
approximation implies that one ignores the presence of overhangs in
the interface. (More often than not, the model is discretized and
$\mathbb R^d, \mathbb R$ are replaced by $\mathbb Z^d,\mathbb Z$.)
Secondly, in the usual spirit of statistical mechanics, the complex
phenomena leading to microscopic interface randomness (e.g. chaotic
motion of water molecules in the case of the ice/water boundary, or the
various atmospheric phenomena determining the motion of individual
snowflakes) are simplified into a probabilistic description where the
dynamics of the height function is modeled by a Markov chain with
simple, ``local'', transition rules.

We already mentioned that on macroscopic scales the interface
evolution looks deterministic: this means that rescaling space as
$\epsilon^{-1}x$, height as $\epsilon h$ and time as
$\epsilon^{-\alpha}t$ (we will discuss the scaling exponent $\alpha>0$
later) and letting $\epsilon\to0$, the random function
$\epsilon h_{\epsilon^{-1}x}(\epsilon^{-\alpha}t)$ converges to a
deterministic function $\phi(x,t)$ that in general is the solution of
a certain non-linear PDE. This is called the \emph{hydrodynamic limit}
and is the analog of a law of large numbers for the sum of independent random
variables. When we say that convergence holds, it does not mean it is
easy to prove it or even to write down the PDE explicitly. Indeed, as
discussed in more detail in the following sections, above dimension
$(d+1)=(1+1)$ the hydrodynamic limit has been proved only for a
handful of models, and one of the goals of this review is to report
on recent results for $d=2$.

On a finer scale, the interface fluctuations around the hydrodynamic limit
are expected to converge, after proper rescaling, to a limit
stochastic process, not necessarily Gaussian. In some situations, but not
always, this is described via a Stochastic PDE. Again, while much is
now known about $(1+1)$-dimensional models (for one-dimensional growth
models and their relation with the so-called KPZ universality class,
we refer to the recent reviews \cite{Corwin,Quastel}), results are
very scarce in higher dimension and we will present some recent ones for $d=2$.

Before entering into more details of the models we consider, it is
important to distinguish between two very different physical
situations. In the case of deposition phenomena, the interface grows
irreversibly and asymmetrically in one direction (say, vertically
upward). The same is the case for the ice/water example if temperature
is not $0$ $^\circ$C: for instance if $T>0$ $^\circ$C then ice melts
and the water phase eventually invades the whole space. In these
situations, the Markov process modeling the phenomenon is irreversible
and the correct scaling for the hydrodynamic limit is the so-called
hyperbolic or Eulerian one: the scaling exponent $\alpha$ introduced above equals $1$. The situation is very
different for the ice/water example exactly at $0$ $^\circ$C: in this
case, the two coexisting phases are at thermal equilibrium and none is
a priori favored. If the ice block occupied a full half-space then the
flat water/ice interface would macroscopically not move and indeed a finite ice
cube evolves only thanks to curvature of its boundary.  In terms of
hydrodynamic limit, one needs to look at longer time-scales than the
Eulerian one: more precisely, one needs to take time of order $\epsilon^{-\alpha},\alpha=2$
(diffusive scaling).

\begin{figure}
\begin{center}
\includegraphics[height=8cm]{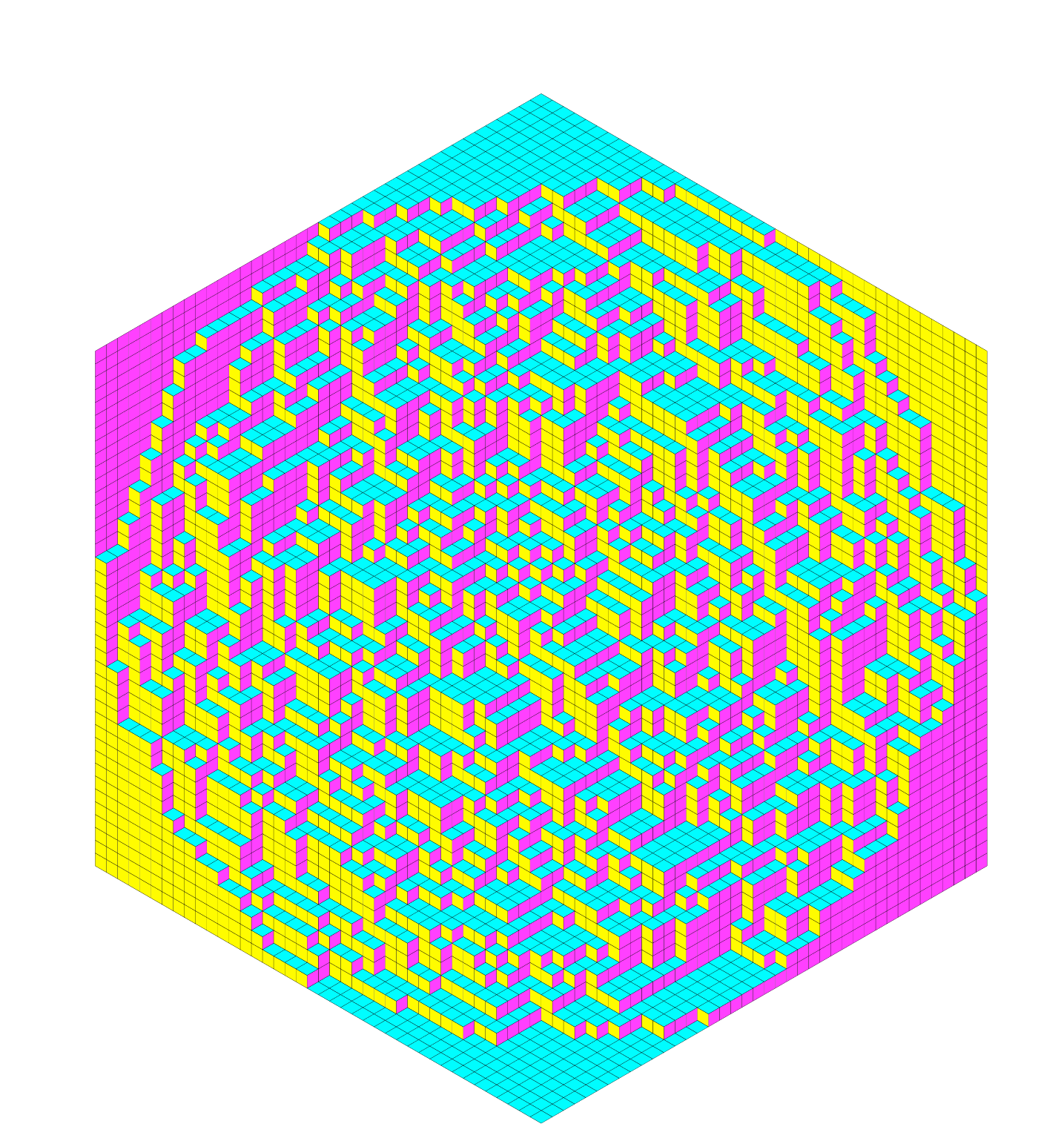}
\caption{A random, uniformly chosen, lozenge tiling of 
a hexagon $\Lambda_\epsilon$  of diameter $\epsilon^{-1}$. As $\epsilon\to0$, the corresponding random interface presents six  ``frozen regions''
  near the corners and a ``liquid region'' inside the arctic circle \cite{Arctic}.}
\label{fig:hexa}
\end{center}
\end{figure}

We will discuss in more detail the Eulerian and diffusive cases,
together with the new results we obtained for some $(2+1)$-dimensional
models, in Sections \ref{sec:introgrowth} and \ref{sec:revers}
respectively. A common feature of all our recent results is that the
interface dynamics we analyze can be formulated as dynamics of dimer
models on bipartite planar graphs, or equivalently of tilings of the plane. See Fig. \ref{fig:hexa} for a 
randomly sampled lozenge tiling of a planar domain.

Such models have a family of
translation-invariant Gibbs measures, with an integrable (actually
determinantal) structure \cite{Kenyon}, that play the role of
stationary states for the dynamics.

\section{Stochastic interface growth}

\label{sec:introgrowth}
\begin{figure}
\begin{center}
\includegraphics[height=3cm]{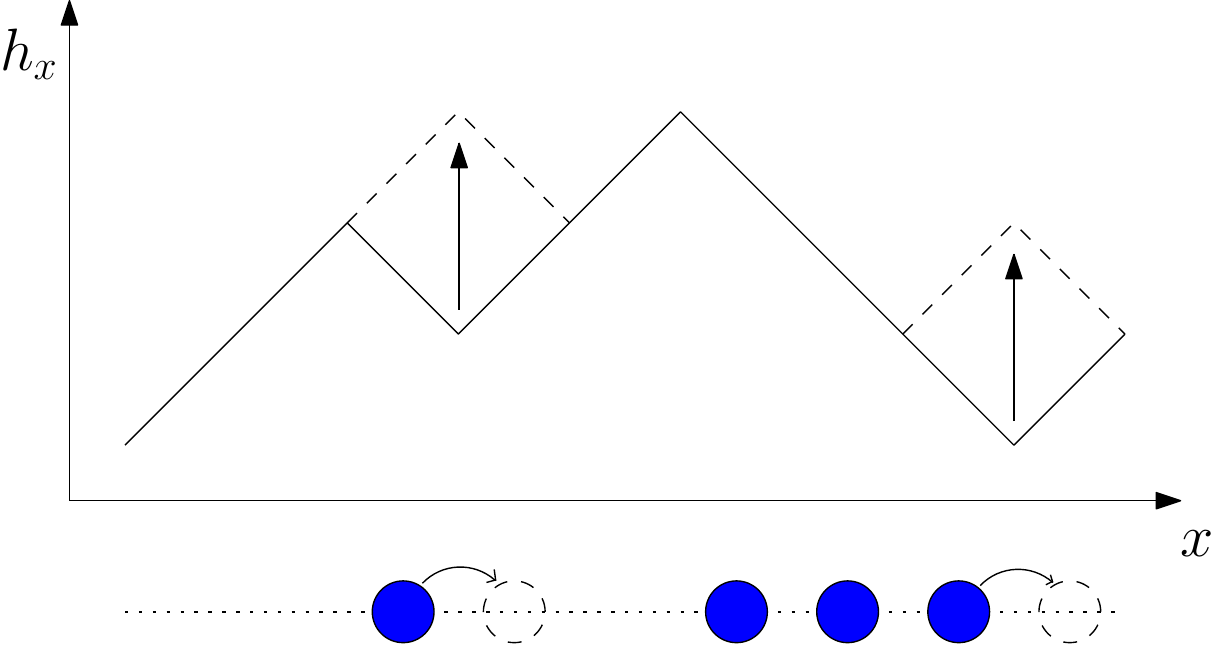}
\caption{In the corner-growth process, heights increase by $1$, with
  transition rate equal to $1$, at local minima. Interpreting a
  negative gradient $h_{x+1}-h_x$ as a particle and a negative one as
  a hole, the dynamics of the gradients is the TASEP: particles try
  independently length-1 jumps to the right, with rate $1$, subject to an
  exclusion constraint (at most one particle per site is allowed). }
\label{fig:TASEP}
\end{center}
\end{figure}

In a stochastic growth process, the height function
$h(t)=\{h_x(t)\}_{x\in \mathbb Z^d}$ evolves asymmetrically, i.e. has
an average non-zero drift, say positive. For instance, growth can be
totally asymmetric: only moves increasing the height are allowed. It
is then obvious that such Markov chain cannot have an invariant
measure. One should look at interface gradients
$\nabla h(t)=\{h_x(t)-h_{x_0}(t)\}_{x\in \mathbb Z^d}$ instead, where
$x_0$ is some reference site (say the origin).  Since the growth
phenomenon we want to model satisfies vertical translation invariance,
the transition rate at which $h_x$ jumps, say, to $h_x+1$ depends only
on the interface gradients (say, the gradients around $x$) and not on
the absolute height $h(t)$. Therefore, the projection of the Markov
chain $h(t)$ obtained by looking at the evolution of $\nabla h(t)$ is
still a Markov chain.  For natural examples one expects that given a
slope $\rho\in \mathbb R^d$, there exists a unique
translation-invariant stationary state $\mu_\rho$ for the gradients,
with the property that $\mu_\rho(h_{x+e_i}-h_x)=\rho_i,i=1,\dots,d$.
A very well known example is the $(1+1)$-dimensional corner growth
model: the evolution of interface gradients is just the 1-dimensional
Totally Asymmetric Simple Exclusion (TASEP), whose invariant measures
are iid Bernoulli product measures labelled by the particle density
$\rho$. See Fig. \ref{fig:TASEP}.

If the initial height profile is sampled from $\mu_\rho$ (more precisely, the height gradients are sampled from $\mu_\rho$, while the height $h_{x_0}$ is assigned some arbitrary value, say zero), then on average the height increases exactly linearly with time:
\begin{eqnarray}
  \label{eq:velocita'}
  \mathbb E_{\mu_\rho}(h_x(t)-h_x(0))=v(\rho)t,\quad v(\rho)>0.
\end{eqnarray}
Now suppose that the initial height profile is instead close to some non-affine profile $\phi_0$, i.e.
\begin{eqnarray}
\label{eq:nowsu}
\epsilon  h_{\epsilon^{-1}x}\stackrel {\epsilon\to0}\to \phi_0(x),\quad \forall x\in \mathbb R.
\end{eqnarray}
Then, one expects that, under so-called hyperbolic rescaling of space-time where $x\to \epsilon^{-1}x,t\to \epsilon^{-1}t$, one has
\begin{eqnarray}
  \label{eq:h1}
 \epsilon h_{\epsilon^{-1} x}(\epsilon^{-1}t)\stackrel[\epsilon\to0]{\mathbb P}{\longrightarrow} \phi(x,t)
\end{eqnarray}
where $\phi(x,t)$ is non-random and solves the 
first order PDE of Hamilton-Jacobi type
\begin{eqnarray}
  \label{eq:h10}
  \partial_t \phi(x,t)=v(\nabla\phi(x,t))
\end{eqnarray}
with $v(\cdot)$ \emph{the same function} as in \eqref{eq:velocita'}.
A couple of remarks are important here:
\begin{itemize}
\item Unless $v(\cdot)$ is a linear function (which is a very
  uninteresting case), the PDE \eqref{eq:h10} develops singularities
  in finite time. Then, one expects $\phi(x,t)$ to solve
  \eqref{eq:h10} in the sense of vanishing viscosity.
  
\item Viscosity solutions of \eqref{eq:h10} are well understood when
  $v(\cdot)$ is convex, since they are given by the variational
  Hopf-Lax formula. However, there is no fundamental reason   why
  $v(\cdot)$ should be convex (this will be an important point in next
  section). Then, much less is known on the analytic side, aside from
  basic properties of existence and uniqueness.
  
\item The example of the TASEP is very special in that invariant
  measures $\mu_\rho$ are explicitly known. One should keep in mind
  that this is an exception rather than the rule and that most
  examples with known stationary measures are $(1+1)$-dimensional. As
  a consequence, the function $v(\cdot)$ in \eqref{eq:h10} is in
  general unknown.
\end{itemize}
Next, let us consider fluctuations in the stationary process started
from $\mu_\rho$. On heuristic grounds, one expects height fluctuations $\hat h_x(t)$
with respect to the average, linear, height profile $\mathbb E_{\mu_\rho}(h_x(t))=\langle\rho, x\rangle+ v(\rho)t$ to
be somehow described, \emph{on large space-time scales}, by a
stochastic PDE (KPZ equation) of the type \cite{KPZ}
\begin{eqnarray}
  \label{eq:SPDE}
  \partial_t \psi(x,t)=\nu\Delta \psi(x,t)+\langle\nabla \psi(x,t),H_{\rho} \nabla \psi(x,t)\rangle+\xi(x,t),
\end{eqnarray}
where:
\begin{itemize}
\item the Laplacian is a diffusion term that tends to locally smooth out fluctuations and $\nu>0$ is a model-dependent constant;
\item the $d\times d$ symmetric matrix $H_\rho$ is the Hessian of the
  function $v(\cdot)$ computed at $\rho$ and
  $\langle\cdot,\cdot\rangle$ denotes scalar product in $\mathbb
  R^d$. This non-linear term  comes just from expanding to second
  order\footnote{The first-order term
    $\langle \nabla v(\rho),\nabla \psi(x,t)\rangle$ in the expansion
    is omitted because it can be absorbed into $\partial_t \psi$ via a
    linear (Galilean) transformation of space-time coordinates.} the
  hydrodynamic PDE \eqref{eq:h10} around the flat solution of slope
  $\rho$.
\item $\xi(x,t)$ is a space-time noise that models the randomness of
  the Markov evolution. It is well known that Eq. \eqref{eq:SPDE} is
  extremely singular if $\xi$ is a space-time white noise (Hairer's
  theory of regularity structures \cite{Hairer} gives a meaning to the
  equation for $d=1$ but not for $d>1$). Since however we are
  interested in properties on large space-time scales and since
  lattice models have a natural ``ultraviolet'' space cut-off of order
  $1$ (the lattice spacing), we can as well imagine that the noise is
  not white in space and its correlation function has instead a decay
  length of order $1$. As a side remark, in the physics literature
  (e.g. \cite{KPZ,Wolf,BarSta}) the presence of a noise regularization
  in space is implicitly understood, and explicitly used in the
  renormalization group computations: this is the cut-off $\Lambda=1$
  that appears e.g. in \cite[App. B]{BarSta}.
  
\end{itemize}
One should not take the above conjecture in the literal sense that the
law of the space-time fluctuation process $\hat h_x(t)$ converges to the law of the solution of
\eqref{eq:SPDE}. Only the large-scale correlation properties of the two should be asymptotically equivalent. 

For $(1+1)$-dimensional models like the TASEP and several others, a
large amount of mathematical results by now supports the following
picture\footnote{As \eqref{eq:SPDE} suggests, for the following to
  hold one needs $v''(\rho)\ne0$, otherwise the fluctuation process
  should be described simply the linear stochastic heat equation with
  additive noise.}: starting say with the deterministic condition
$h_x(0)\equiv0$, the standard deviation of $h_x(t)$ grows as
$t^{\beta},\beta=1/3$, the space correlation length grows like
$t^{1/z}$, where $z=3/2$ is the so-called dynamic exponent and the fluctuation field $\hat h_x(t)$ rescaled
accordingly tends as $t\to\infty$ to a (non-Gaussian) limit
process. We do not enter into any more detail for $(1+1)$-dimensional
models of the KPZ class here, see for instance the reviews
\cite{Quastel,Corwin}; let us however note that this behavior is very
different from the (Gaussian) one of the stochastic heat equation with
additive noise (called ``Edwards-Wilkinson equation'' in the physics
literature), obtained by dropping the non-linear term in
\eqref{eq:SPDE}.

On the other hand, for $(d+1)$-dimensional models, $d\ge 3$,
renormalization-group computations \cite{KPZ} applied to the
stochastic PDE \eqref{eq:SPDE} suggest that, if the non-linear term is
sufficiently small (in terms of the microscopic growth model: if the
speed function $v(\cdot)$ is sufficiently close to an affine function)
then non-linearity is irrelevant, meaning that the large-scale
fluctuation properties of the model (or of the solution of
\eqref{eq:SPDE}) are asymptotically the same as those of the
stochastic heat equation: these models belong to the so-called
Edwards-Wilkinson universality class. There is very recent
mathematical progress in this direction: indeed, \cite{MU} states that
for $d\ge 3$ the solution of \eqref{eq:SPDE} tends on large space-time
scales to the solution of the Edwards-Wilkinson equation, if
$H_\rho=\lambda \mathbb I$ with $\mathbb I$ the $d\times d$ identity
matrix and $\lambda$ small enough. See also \cite{Z} where similar
results are stated for the $d\ge3$ dimensional stochastic heat
equation with multiplicative noise, that is obtained from
\eqref{eq:SPDE} via the Cole-Hopf transform.
  
The situation is richer in the borderline case of the critical dimension $d=2$, to which the next two sections are devoted.

\subsection{$(2+1)$-dimensional growth: KPZ and  Anisotropic KPZ (AKPZ) classes}
For $(1+1)$-dimensional models, the non-linear term in \eqref{eq:SPDE}
equals $v''(\rho)(\partial_x \psi(x,t))^2$: multiplying $\psi$ by a
suitable  constant, we can always
replace $v''(\rho)\ne0$ by a positive  constant. The picture is richer
for $d>1$, and in particular in the case $d=2$ we consider here. In
fact, one should distinguish two cases:
\begin{enumerate}
\item (Isotropic) KPZ class:  $\det(H_\rho)>0$   (strictly);
  
\item Anisotropic KPZ (AKPZ) class:  $\det(H_\rho)\le0$.
\end{enumerate}
According to whether a growth model has a speed function $v(\cdot)$
whose Hessian satisfies the former or latter condition, the
large-scale behavior of its fluctuations is conjectured to be very
different.

The {\bf isotropic KPZ class} is the one considered in the original
KPZ work \cite{KPZ}. In this case, perturbative renormalization-group
arguments\footnote{``perturbative'' here means that, if we imagine that
  the non-linear term in \eqref{eq:SPDE} has a prefactor $\lambda$,
  then one expands the solution around the linear $\lambda=0$ case,
  keeping only terms up to order
  $O(\lambda^2)$
  .}  suggest that fluctuations of $h_x(t)$ (or of the solution
$\psi(x,t)$ of \eqref{eq:SPDE}) grow in time like $t^\beta$ and that,
in the stationary states, fluctuations grow in space as
${\rm Var}_{\mu_\rho} (h_x-h_y)\sim |x-y|^{2\alpha}$, with two
exponents $\beta>0,\alpha=2\beta/(\beta+1)$ that are \emph{different}
from those of the Edwards-Wilkinson equation: non-linearity is said to
be \emph{relevant}\footnote{ The relation $\alpha=2\beta/(\beta+1)$
  is another way of writing a scaling relation between exponents that
  is usually written as $\alpha+z=2$ where $z=\alpha/\beta$. Here $z$
  is the so-called dynamic exponent that equals $3/2$ for
  one-dimensional KPZ models.}. The Edwards-Wilkinson equation can be
solved explicitly and in two dimensions one finds $\alpha_{EW}=\beta_{EW}=0$
(growth in time and space is only logarithmic; the stationary state is
the (log-correlated) massless Gaussian field). The values of
$\alpha,\beta$ for the isotropic KPZ equation cannot be guessed by
perturbative renormalization-group arguments and they are accessible
only through numerical experiments (see discussion below). Note that
$\alpha>0$ means that stationary height profiles are much rougher
than a lattice massless Gaussian field.

The {\bf Anisotropic KPZ} case was analyzed later by Wolf \cite{Wolf}
with the same renormalization-group approach and the result came out as
a surprise: non-linearity turns out to be non-relevant in this case, i.e., the
growth exponents $\alpha,\beta$ are predicted to be $0$ as for
the Edwards-Wilkinson equation.

Let us summarize this discussion into
a conjecture:
\begin{conj}
  \label{conj:Wolf}
  Let $v(\cdot)$ be the speed function of a (reasonable)
  $(2+1)$-dimensional growth model. If $\det(H_\rho)>0$ with $H_\rho$
  the Hessian of $v(\cdot)$ computed at $\rho$, then height fluctuations grow in time as
  $t^\beta$ for some model-independent $\beta>0$ and height
  fluctuations in the stationary states grow as distance to the power
  $2\beta/(\beta+1)$.  If instead $\det(H_\rho)\le0$, then
  $\beta=\alpha=0$ and the stationary states have the same height correlations in space as a massless Gaussian field.
\end{conj}

Let us review the evidence in favor of this conjecture, apart from the
renormalization-group argument of \cite{KPZ,Wolf} that does not provide much intuition and seems very hard
to be turned into a mathematical proof:
\begin{enumerate}

\item A somewhat rough but suggestive argument that sheds some light
  on Conjecture \ref{conj:Wolf} is given in \cite[Sec. 2.2]{Prahofer}. One
  imagines that in the evolution of the fluctuation field $\psi$ there are two
  effects. Thermal noise adds random positive or negative ``bumps'',
  at random times, to the initially flat height profile; each bump then evolve following the hydrodynamic
  equation, expanded to second order:
  \begin{eqnarray}
    \label{eq:bumps}
  \partial_t\psi=\langle\nabla \psi(x,t),H_{\rho} \nabla
  \psi(x,t)\rangle.   
  \end{eqnarray}
  It is not hard to convince oneself that, under \eqref{eq:bumps}, if
  both eigenvalues of $H_\rho$ are, say, strictly positive, then a
  positive bump grows larger with time and a negative bump shrinks (the reverse happens if the eigenvalues of $H_\rho$ are both negative).
  See Figure \ref{fig:bump}. On the other hand, if $\det(H_\rho)<0$ (so that the eigenvalues of
  $H_\rho$ have opposite signs)  then a positive bump spreads
  in the direction where the curvature of $v(\cdot)$ is positive, but
  its height shrinks because of the concavity of $v(\cdot)$ in the
  other direction (the same argument applies to negative bumps). Then,
  it is intuitive that when $\det(H_\rho)<0$ height fluctuations
  should grow slower with time than when $\det(H_\rho)>0$, where the
  effects of spreading positive bumps accumulate.

\begin{figure}
\begin{center}
\includegraphics[height=4cm]{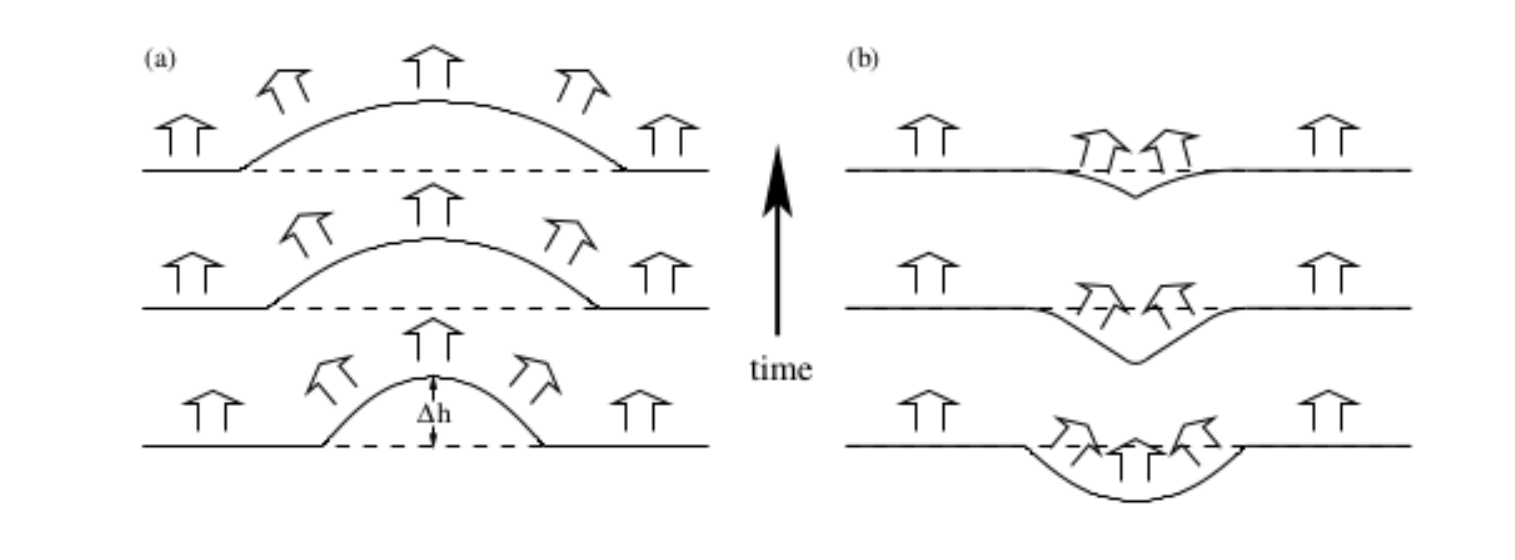}
\caption{Left: The time evolution of a positive bump   under equation \eqref{eq:bumps}, when $\det(H_\rho)>0$. 
The height $\Delta h$ of the bump is constant in time  while its width grows as $t^{1/2}$. Right:  a negative bump, on the other hand,  develops a cusp and its height decreases as $1/t$. This figure is taken from \cite{Prahofer}. }
\label{fig:bump}
\end{center}
\end{figure}

\item There exist some growth models that satisfy a so-called
  ``envelope property'', saying essentially that given two initial
  height profiles $\{h^{(j)}_x\}_{x\in\mathbb Z^d}, j=1,2$,  one can find a coupling between  the
  corresponding profiles $\{h^{(j)}_x(t)\}_{x\in\mathbb Z^d}$ at time
  $t$ such that the evolution started from
  the profile
  $\check h:=\{\max(h^{(1)}_x,h^{(2)}_x)\}_{x\in\mathbb Z^d}$ equals
  $\check h(t)=\{\max(h^{(1)}_x(t),h^{(2)}_x(t))\}_{x\in\mathbb
    Z^d}$. One example is the $(2+1)$-dimensional corner-growth model
  analogous to that of Fig. \ref{fig:TASEP} except that the interface is two-dimensional and unit cubes
  instead of unit squares are deposed with rate one on it.
  For growth models satisfying the envelope property, a
  super-additivity argument implies that the hydrodynamic limit
  \eqref{eq:h1} holds and moreover that the function $v(\cdot)$ in
  \eqref{eq:h10} is convex \cite{Seppa,Reza}. While for
  $(2+1)$-dimensional models in this class the stationary measures
  $\mu_\rho$ and the function $v(\cdot)$ cannot be identified
  explicitly, convexity implies  (at least in the region of slopes 
  where $v(\cdot)$ is smooth and strictly convex) that 
  $\det(H_\rho)>0$: these models must belong to the isotropic KPZ
  class. The $(2+1)$-dimensional corner-growth model was studied
  numerically in \cite{Tang} and it was found, in agreement with
  Conjecture \ref{conj:Wolf}, that $\beta\simeq 0.24$ (the numerics is
  sufficiently precise to rule out the value $1/4$ which was
  conjectured in earlier works). The same value for $\beta$ is found
  numerically \cite{HHA92} from direct simulation of (a space
  discretization of) the stochastic PDE \eqref{eq:SPDE} with $\det(H_\rho)>0$.
  
\item For models in the AKPZ class there is no chance to get the
  hydrodynamic limit by simple super-additivity arguments since, as we
  mentioned, $v(\cdot)$ would turn out to be convex. On the other
  hand, as we discuss in more detail in next section, there exist some
  $(2+1)$-dimensional growth models for which the stationary measures
  $\mu_\rho$ can be exhibited explicitly, and they turn out to be of
  massless Gaussian type, with logarithmic growth of fluctuations:
  $\alpha=0$. For such models, one can prove also that $\beta=0$ and
  one can compute the speed function $v(\cdot)$. In all the known
  examples, a direct computation shows that $\det(H_\rho)<0$, as it should according to Conjecture \ref{conj:Wolf}.

\end{enumerate}

\begin{rem}
\label{rem:emph}
  Let us emphasize that, in general, it is not possible to read a
  priori, from the generator of the process, the convexity properties
  of the speed function $v(\cdot)$, and therefore its universality
  class. This is somehow in contrast with the situation in equilibrium
  statistical mechanics, where usually the universality class of a
  model can be guessed from symmetries of its Hamiltonian.  It is even conceivable, though we are not aware of any concrete example, that there exist
  growth models for which the sign of $\det(H_\rho)$ depends on $\rho$.
\end{rem}

\subsection{Mathematical  results for Anisotropic KPZ growth models}

\label{sec:newAKPZ}
As we already mentioned, there are no results other than numerical
simulations or non-rigorous arguments supporting the part of
Conjecture \ref{conj:Wolf} concerning the isotropic KPZ
class. Fortunately, the situation is much better for the AKPZ class,
which includes several models that are to some extent ``exactly
solvable''.

\begin{figure}
\begin{center}
\includegraphics[height=5cm]{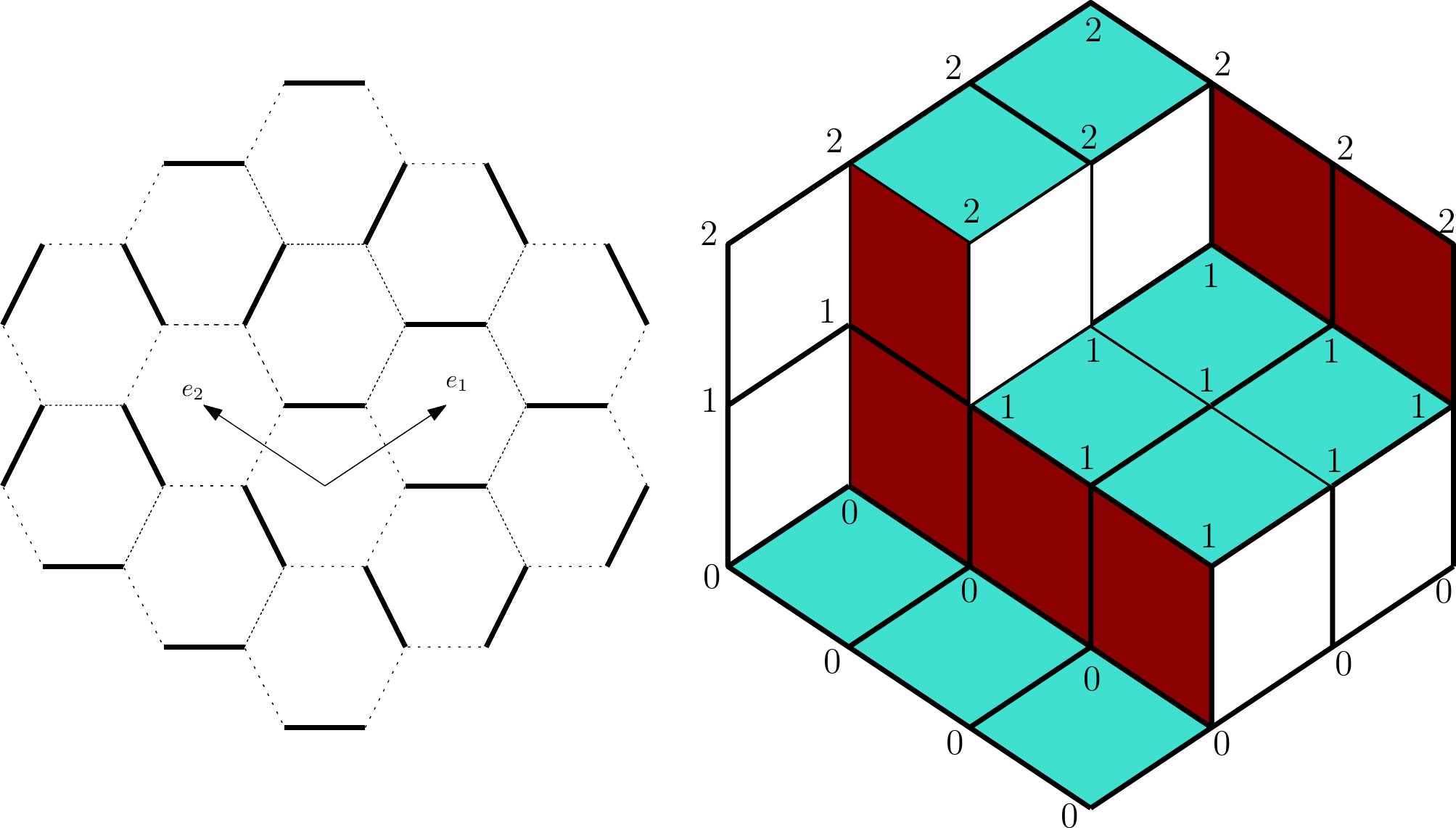}
\caption{A perfect matching of the hexagonal lattice (left) and the
  corresponding lozenge tiling (right). Near each lozenge vertex is
  given the height of the interface w.r.t. the horizontal plane. For
  clarity let us emphasize that, while we draw only a finite portion
  of the matching/tiling, one should imagine that it extends to a
  matching/tiling of the infinite graph/plane.}
\label{fig:mapping}
\end{center}
\end{figure}

Several of the AKPZ models for which mathematical results are
available have a height function that can be associated to a
two-dimensional dimer model (an exception is the Gates-Westcott model
solved by Pr\"ahofer and Spohn \cite{PS97}).  Let us briefly recall
here a few well-known facts on dimer models (we refer to \cite{Kenyon}
for an introduction). For definiteness, we will restrict our
discussion to the dimer model on the infinite hexagonal graph but most
of what we say about the height function and translation-invariant
Gibbs states extends to periodic, two-dimensional bipartite graphs
(say, $\mathbb Z^2$).  A (fully packed) dimer configuration is a
perfect matching of the graph, i.e., a subset $M$ of edges such that each vertex of the graph is contained in one and only one edge in $M$; as in Fig. \ref{fig:mapping}, in the
case of the hexagonal graph the matching can be equivalently seen as a
lozenge tiling of the plane and also as a monotone discrete
two-dimensional interface in three dimensional space.  ``Monotone''
here means that the interface projects bijectively on the plane
$x+y+z=0$. The height function is naturally associated to vertices of
lozenges, i.e. to hexagonal faces. We will use the dimer, the tiling or the
height function viewpoint interchangeably.
If on the graph we choose coordinates $x=(x_1,x_2)$ according to the
axes $e_1,e_2$ drawn in Fig. \ref{fig:mapping}, it is easy to see that
the overall slope $\rho=(\rho_1,\rho_2)$ of the interface must belong
to the triangle $ \mathbb T\subset \mathbb R^2$ defined by the inequalities
$0\le \rho_1\le 1,0\le \rho_2\le 1, 0\le \rho_1+\rho_2\le 1$. It is
known that, given $\rho$ in the interior of $\mathbb T$, there exists
a unique translation-invariant ergodic Gibbs state $\pi_\rho$ of slope
$\rho$. That is, $\pi_\rho$ is a (translation invariant, ergodic)
probability measure on dimer configurations of the infinite graph,
such that the average height slope is $\rho$ and such that,
conditionally on the configuration outside any finite domain
$\Lambda$, the law of the configuration inside $\Lambda$ is uniform
over all dimer configurations compatible with the outside (DLR
condition). In fact, much more is known: as a consequence of
Kasteleyn's theory \cite{Kasteleyn}, such measures have a determinantal
representation. That is, the probability of a cylindrical event of the
type ``$k$ given edges are occupied by dimers'' is given by the
determinant of a $k\times k$ matrix, whose elements are the Fourier
coefficients of an explicit function on the two-dimensional torus
$\{(z,w)\in \mathbb C^2: |z|=|w|=1\}$ \cite{KOS}.  Thanks to this representation,
much can be said about large-scale properties of the measures
$\pi_\rho$. Notably, correlations decay like the inverse distance
squared and the height function  scales to a massless Gaussian
field with logarithmic covariance structure.

\begin{figure}
\begin{center}
\includegraphics[height=5cm]{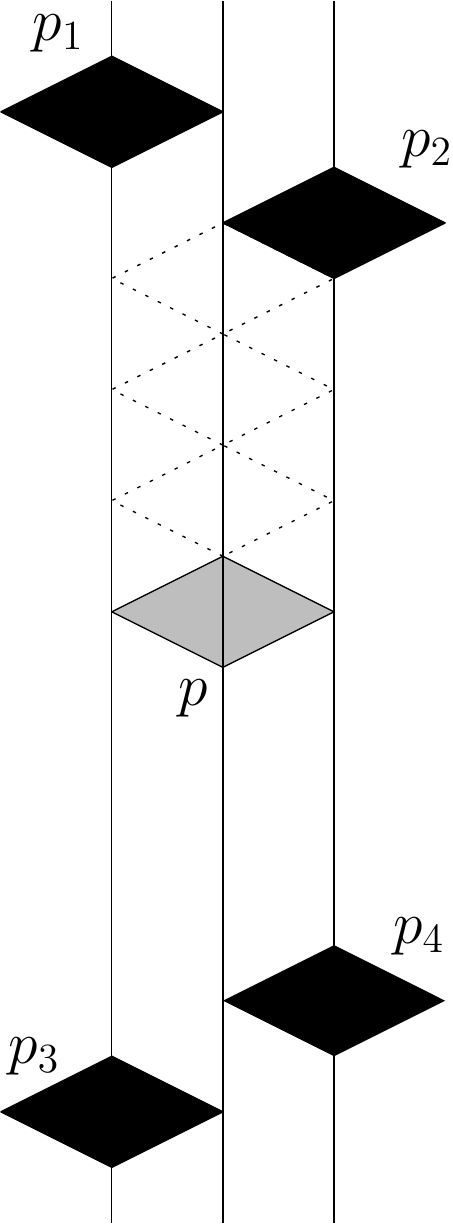}
\caption{Each particle (or horizontal lozenge) $p$ is constrained between its four neighboring particles $p_1,\dots p_4$. The three positions particle $p$ can jump to (with rate $1$) in the Borodin-Ferrari dynamics are dotted.}
\label{fig:interlac}
\end{center}
\end{figure}

Now that we have a nice candidate for a $(2+1)$-dimensional height
function, we go back to the problem of defining a growth model that
would hopefully be mathematically treatable and shed some light on Conjecture
\ref{conj:Wolf}. To this purpose, let us remark first of all that, to
a lozenge tiling as in Fig. \ref{fig:mapping}, one can bijectively
associate a \emph{two-dimensional system of interlaced particles}. For
this purpose, we will call ``particles'' the horizontal (or blue)
lozenges (the positions of the others are uniquely determined by
these) and we note that particle positions along a vertical column are
interlaced with those of the two neighboring columns. See
Fig. \ref{fig:interlac}.
A first natural candidate for a growth process would be the following
immediate generalization of the TASEP: each particle jumps $+1$
vertically, with rate $1$, provided the move does not violate the
interlacement constraints. Actually, this is nothing but the
three-dimensional corner-growth model. As we already mentioned, this
should belong to the isotropic KPZ class and its stationary measures
$\mu_\rho$ should be extremely different from the Gibbs measures
$\pi_\rho$, with a power-like instead of logarithmic growth of
fluctuations in space. Unfortunately, none of this could be
mathematically proved so far.

In the work \cite{BF08}, A. Borodin and P. Ferrari
considered instead another totally asymmetric growth model where
each particle can jump an unbounded distance $n$ upwards, with rate
independent of $n$ (say, rate $1$), provided the interlacements are still satisfied
after the move. See Fig. \ref{fig:interlac}. The situation is then
entirely different with respect to the corner-growth process: the two processes belong to two different universality classes. If the
initial condition of the process is a suitably chosen, deterministic, fully packed
particle arrangement (see Fig. 1.1 in \cite{BF08}), it was shown that
the height profile rescaled as in \eqref{eq:h1} does converge to a
deterministic limit $\phi(\cdot,t)$, that solves the Hamilton-Jacobi equation
\eqref{eq:h10} with
\begin{eqnarray}
  \label{eq:vmu}
  v(\rho)=\frac1\pi\frac{\sin(\pi \rho_1)\sin(\pi\rho_2)}{\sin(\pi(\rho_1+\rho_2))}.
\end{eqnarray}
A couple of remarks are important for the subsequent discussion:
\begin{itemize}
\item with the initial condition chosen in \cite{BF08}, $\phi(x,t)$ turns out to be a \emph{classical}  solution of \eqref{eq:h10}. That is, the characteristic lines of the PDE do not cross at positive times (we emphasize that this is due
to the specific form of the chosen initial profile and not to the form of $v(\cdot)$).
\item An explicit computation shows that $\det(H_\rho)<0$ for the
  function \eqref{eq:vmu}: this growth model is then a candidate to
  belong to the AKPZ class. 
\end{itemize}
As mentioned in Remark \ref{rem:emph} above, let us emphasize that we see no obvious way to guess a priori that the corner growth process and the ``long-jump'' one should belong to different universality classes.

Various other results were proven in \cite{BF08}, but let us mention
only two of them, that support the conjecture that this model indeed belongs
to the AKPZ class:
\begin{enumerate}
\item the fluctuations of
$h_{\epsilon^{-1}x}(\epsilon^{-1}t)$ around its average value are of order
$\sqrt{\log 1/\epsilon}$ (the growth exponent is $\beta=0$) and, once rescaled by this factor, they tend
to a Gaussian random variable;
\item the local law of the interface
gradients at time $\epsilon^{-1}t$ around the point $\epsilon^{-1}x$
tends to the Gibbs measure $\pi_\rho$ with $\rho=\nabla\phi(x,t)$. 
\end{enumerate}
\begin{rem}
\label{rem:special}
  The basic fact behind the results of \cite{BF08} is that for the
  specific choice of initial condition, one can write \cite[Th. 1.1]{BF08} the probability
  of certain events of the type ``there is a particle at position
  $x_i$ at time $t_i, i\le k$'' as a $k\times k$ determinant, to which
  asymptotic analysis can be applied. The same determinantal
  properties hold for other ``integrable'' initial conditions, but
  they are not at all a generic fact.
\end{rem}
Point (2)  above clearly suggests that the Gibbs measures $\pi_\rho$ should
be  stationary states for the interface gradients. In fact, this
is a result I later proved:
 \begin{thm}\cite[Th. 2.4]{T2+1}
   For every slope $\rho$ in the interior of $\mathbb T$, the measure
   $\pi_\rho$ is stationary for the process of the interface
   gradients.
 \end{thm}
 Recall that,
as discussed above, the Gibbs measures $\pi_\rho$ of the dimer model
have the large-scale correlation structure of a massless Gaussian
field and indeed Conjecture \ref{conj:Wolf} predicts that stationary
states of AKPZ growth processes behave like massless fields.  Most of
the technical work in \cite{T2+1} is related to the fact that, since
particles can perform arbitrarily long jumps with a rate that
\emph{does not decay} with the jump length, it is not clear a priori
that the process exists at all: one can exhibit initial configurations
such that particles jump to $+\infty$ in finite time (this issue does
not arise in the work \cite{BF08} where, thanks to the chosen initial
condition, there is no difficulty in defining the infinite-volume
process). In \cite{T2+1} it is shown via a comparison with the
one-dimensional Hammersley process \cite{SeppaH} that, for a typical
initial condition sampled from $\pi_\rho$, particles jump almost
surely a finite distance in finite time and that, despite the
unbounded jumps, perturbations do not spread instantaneously through
the system. This means that if two initial configurations differ only on a subset $S$ of the lattice, their evolutions can be coupled so that at finite time $t$ they are with high probability equal sufficiently far away from $S$ (how far, depending on $t$).

To follow the general program outlined above, once the stationary
states are known, one would like to understand the growth exponent
$\beta$ for the stationary process. We proved the following, implying $\beta=0$:
\begin{thm}
 \cite[Th. 3.1]{T2+1}  For every lattice site $x$, we have
 \begin{eqnarray}
\label{eq:u}
   \limsup_{t\to\infty}\mathbb P_{\pi_\rho}\left(|h_x(t)-\mathbb E_{\pi_\rho}(h_x(t))|\ge u\sqrt{\log t}\right)\stackrel{u\to\infty}\longrightarrow 0
 \end{eqnarray}
where
\begin{eqnarray}
  \label{eq:a}
\mathbb E_{\pi_\rho}(h_x(t))=v(\rho)t+\langle x,\rho\rangle.  
\end{eqnarray}
\end{thm}
To be precise:
\begin{itemize}
\item in the statement of
\cite[Th. 3.1]{T2+1} there is a technical restriction on the slope
$\rho$, that was later removed in joint work with S. Chhita and
P. Ferrari \cite{CFT};
\item the proof that the speed of
  growth $v(\cdot)$ in \eqref{eq:a} is the same as the function
  $v(\cdot)$ in \eqref{eq:vmu}, as it should, was obtained by
  S. Chhita and P. Ferrari in \cite{CF15} and requires a nice
  combinatorial property of the Gibbs measures $\pi_\rho$.
\end{itemize}
With reference to Remark \ref{rem:special} above, it is important to
emphasize that there is no known determinantal form for the space-time
correlations of the stationary process; for the proof of
\eqref{eq:u} we used a more direct and probabilistic method.


Finally, it is natural to try to obtain a hydrodynamic limit for the
height profile. Recall that in \cite{BF08} such a result was proven
for an ``integrable'' initial condition that allowed to write certain
space-time correlations, and as a consequence the  average particle
currents, in determinantal form. On the other hand, convergence to the
hydrodynamic limit should be a very robust fact and not rely on such special structure. We have indeed:
\begin{thm}\cite[Th. 3.5 and 3.6]{LegrasT}
Let the initial height profile satisfy \eqref{eq:nowsu}, with $\phi_0$ a Lipshitz function with gradient in the interior of $\mathbb T$. 
Let one of the following two conditions be satisfied:
\begin{itemize}
\item $\phi_0$ is $C^2$ and the time $t$ is smaller than $T$, the maximal time up to which a classical solution of \eqref{eq:h10} exists;
\item $\phi_0$ is either convex or concave (in which case we put no restriction on $t$). 
\end{itemize}
Then, the convergence \eqref{eq:h1} holds, with $\phi(x,t)$ the viscosity solution of \eqref{eq:h10}.
\end{thm}
The restriction to either small times or to convex/concave profile is due to the
fact that we have in general little analytic control on the singularities of
\eqref{eq:h10}, due to the non-convexity of $v(\cdot)$.  For convex
initial profile,  the viscosity solution of the PDE is given by a Hopf variational form 
and this allows to bypass these analytic difficulties.  Let us emphasize, to
avoid any confusion, that even in the case of convex initial profile
the solution \emph{does} in general develop singularities (shocks),
i.e. discontinuities in space of the gradient $\nabla \phi(x,t)$.
\begin{op}
  Are  height fluctuations still $O(\sqrt{\log t})$ at the location of shocks?
\end{op}
Another important observation is the following. Given that we know
explicitly the stationary states of the process and that the dynamics
is monotone (i.e, if an initial profile is higher than another, under a suitable coupling it
will stay higher as time goes on), it is tempting to try to apply the
method developed by Rezakhanlou in \cite{Reza2}, that gives under such
circumstances convergence of the height profile of a growth model to the
viscosity solution of the limit PDE.  The delicate point is however
that \cite{Reza2} crucially requires that perturbations spread at
finite speed through the system, so that one can analyze the evolution
``locally'', in small enough windows where the profile can be
approximated by  one sampled from $\mu_\rho$, with suitably chosen
slope $\rho$ that depends on the window location.  Due to
unboundedness of particle jumps, however, the ``finite-speed propagation
property'' might fail in our case and in any case it cannot hold
uniformly for all initial conditions.  Most of the technical work in
\cite{LegrasT} is indeed devoted to proving that one can localize the
dynamics despite the long jumps.  A crucial fact is that we show that
the growth process under consideration can be reformulated through a
so-called Harris-like graphical construction.


\subsubsection{Extensions and open problems}
\label{sec:ext}
\begin{figure}
\begin{center}
\includegraphics[height=5cm]{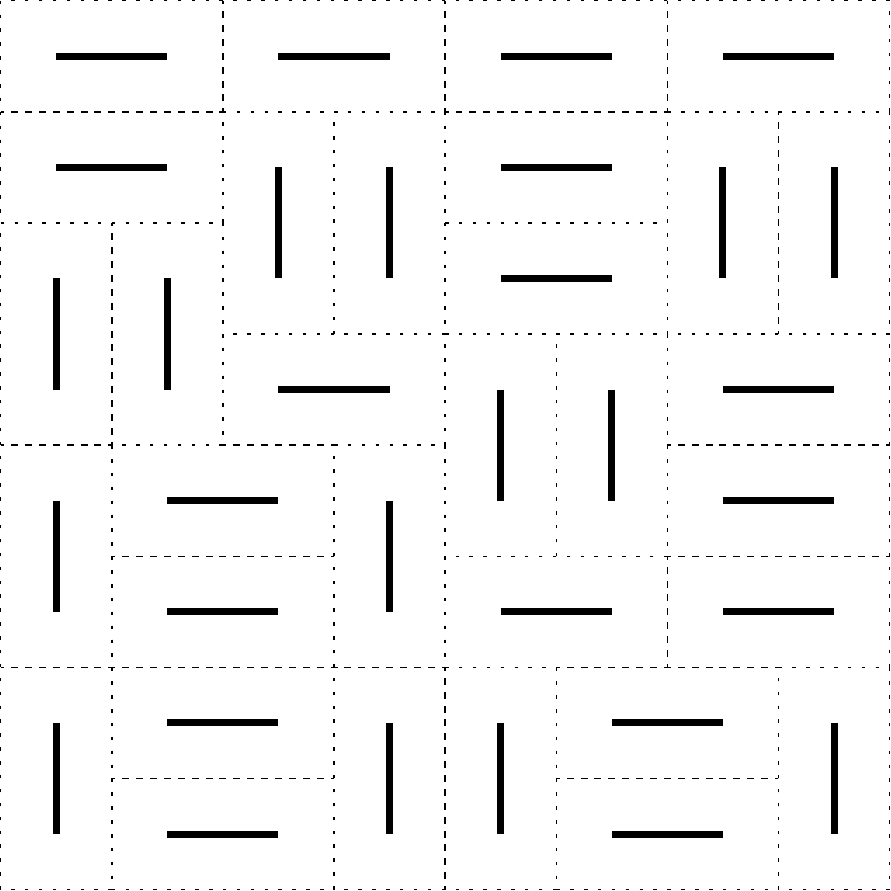}
\caption{A perfect matching of the square lattice and the corresponding domino tiling (dotted). See \cite[Fig. 2]{CFT} for the definition of ``particles'' and their interlacing relations, and \cite{Kenyon} for the definition of the height function.}
\label{fig:domino}
\end{center}
\end{figure}

There are various ways how the ``lozenge tiling dynamics with long
particle jumps'' of previous section can be generalized to provide
other $(2+1)$-dimensional growth processes in the AKPZ class. One such
generalization was given in \cite[Sec. 3.1]{T2+1}. There, one starts
with the observation that: (i) as was the case for lozenge tilings,
also domino tilings of $\mathbb Z^2$ (dominoes being $2\times 1$
rectangles, horizontal or vertical, see Fig. \ref{fig:domino}) have a natural height function
interpretation, and (ii) a domino tiling can be bijectively mapped to
a two-dimensional system of interlaced particles (interlacement
constraints are different than in lozenge case).  

This suggests a
growth process where particles jump in an asymmetric fashion and the
transition rate is independent of the jump length, jumps being limited
only by the interlacement constraints. Then, the same results that
were proven for the lozenge dynamics (notably, stationarity of the
Gibbs measures, logarithmic correlations in space in the stationary
states ($\alpha=0$) and logarithmic growth of fluctuations in the
stationary process implying $\beta=0$) hold in this case too. The
speed of growth $v(\cdot)$ for the domino dynamics was later computed
in a joint work with S. Chhita and P. Ferrari \cite[Th. 2.3]{CFT}: it
turns out to be rather more complicated than \eqref{eq:vmu}, but it is
still an explicit function for which one can prove with some effort
that the Hessian has negative determinant, in agreement with Conjecture
\ref{conj:Wolf}.

Finally, there is yet another class of driven two-dimensional
interlaced particle systems, that was introduced in \cite{BBO}. While
these have rather a group-theoretic motivation, these processes can
also be viewed as $(2+1)$-dimensional growth models and actually the
main result of \cite{BBO} can be seen as a hydrodynamic limit for the
height function \cite[Sec. 3.3]{BBO}.  Once again, direct inspection
of the Hessian of the velocity function shows that these models belong
to the AKPZ class\footnote{For the growth models of \cite{BBO}, the
  determinant of the Hessian of the speed was computed by Weixin Chen,
  as mentioned in the unpublished work
  http://math.mit.edu/research/undergraduate/spur/documents/2012Chen.pdf
}, so these provide other natural candidates where Wolf's prediction
of logarithmic growth of fluctuations can be tested (the logarithmic
nature of fluctuation correlations is conjectured in \cite{BBO}; we
are not aware of an actual proof).

In conclusion, there are now quite a few $(2+1)$-dimensional growth
models in the AKPZ class for which Wolf's predictions in Conjecture
\ref{conj:Wolf} can be verified. There is however one aspect one may find
rather unsatisfactory.  Both for the lozenge tiling dynamics, where
the speed function turns out to be given by \eqref{eq:vmu} and for its
domino tiling generalization, where $v(\cdot)$ is a much more complicated-looking
combination of ratios of trigonometric functions (see
\cite[Eq. (2.6)]{CFT}) and also for the interlaced particle dynamics
of \cite{BBO}, one verifies via brute-force computation that the
Hessian of $H_\rho$ of the corresponding velocity function $v(\cdot)$
has negative determinant. The frustrating fact is that via the explicit
computation one does not see at all how the sign of the determinant the Hessian is
related  to the model being in the Edwards-Wilkinson
universality class! We are still far from having a meta-theorem
saying ``if the exponents $\alpha$ and $\beta$ are zero, then the
determinant of the Hessian is negative''. Up to now, we have essentially
heuristic arguments and ``empirical evidence'' based on a few mathematically treatable models.

\begin{op} It would be very interesting to prove that the Hessian of
  the velocity function for the growth models just mentioned has
  negative determinant \emph{without} going through the explicit
  computation of the second derivatives.
\end{op}

\subsubsection{Slow decorrelation along the characteristics}

The results we discussed in Sections \ref{sec:newAKPZ} and
\ref{sec:ext} (growth of fluctuation variance with time and spatial
correlations in the stationary state) concern fluctuation properties
\emph{at a single time}. Another question of great interest is how
fluctuations at different space-time points $(x_i,t_i)$ are
correlated. For $(1+1)$-dimensional growth models in the KPZ class,
the following picture has emerged \cite{FerrariSlow}: correlation
decay slowly along the characteristic lines of the PDE \eqref{eq:h10}, and faster
along any other direction.  For instance, take two space-time points
$(x_1,t_1)$ and $(x_2,t_2), t_1<t_2$ and think of $t_2$ large. If the
two points are on the same characteristic line, then the height
fluctuations (divided by the rescaling factor $t_i^{\beta}=t_i^{1/3}$)
will be almost perfectly correlated as long as $t_2-t_1\ll t_2$. If
instead the two points are \emph{not} along a characteristic line,
then correlation will be essentially zero as soon as
$t_2-t_1\gg t_2^{1/z}, z=3/2$. 

It has been conjectured \cite{FerrariSlow,BF08} that a similar
phenomenon of slow decorrelation along the characteristic lines should
occur for $(2+1)$-dimensional growth.  For the AKPZ models described
in the previous sections, it is still an open problem to prove
anything in this direction.  In the work \cite{BCT} in collaboration
with A. Borodin and I. Corwin, we studied a growth model that depends
on a parameter $q\in[0,1)$: for $q=0$ it reduces to the long-jump
lozenge dynamics of \cite{BF08,T2+1}, while if $q\to1$ and particle
distances are suitably rescaled the dynamics simplifies in that
fluctuations become Gaussian. In this limit, we were able to prove
that, if height fluctuations are computed along characteristic lines,
their correlations converge to those of the Edwards-Wilkinson equation
and in particular they are large as long as $t_2-t_1\ll t_2$. If
correlations are computed instead along a different direction, then
they are essentially zero as soon as $t_2-t_1\gg t_2^{1/z}$, where
$z=2$ is the dynamic exponent of the Edwards-Wilkinson equation.

%




\section{Interface dynamics at thermal equilibrium}
 \label{sec:revers}
 Let us now move to reversible interface dynamics (we refer to
 \cite{SpohnJSP,Funaki} for an introduction). We can
 imagine that the interface is defined on a finite subset
 $\Lambda_\epsilon$ of $\mathbb Z^d$ of diameter $O(\epsilon^{-1})$,
 say the cubic box $[0,\dots,\epsilon^{-1}]^d$ so that after the
 rescaling $x=\epsilon^{-1}\xi$, the space coordinate $\xi$ is in the
 unit cube. We impose Dirichlet boundary conditions, i.e., for
 $x\in \partial \Lambda_\epsilon$ the height $h_x(t)$ is fixed to some
 time-independent value $\bar h_x$. The way to model the evolution of
 a phase boundary at thermal equilibrium  is
 to take a Markov process with stationary and reversible measure of
 the Boltzmann-Gibbs form (we absorb the inverse temperature into the potential $V$)
\begin{eqnarray}
  \label{eq:BG}
  \pi_{\Lambda_\epsilon}(h)\propto e^{-1/2\sum_{x\sim y} V(h_x-h_y)}
\end{eqnarray}
where the sum runs, say,  on nearest neighboring pairs of  vertices. Note that the
potential $V$ depends only on interface gradients and not on the
absolute height itself: this reflects the vertical translation
invariance of the problem (apart from  boundary
condition effects). A minimal requirement on $V$ is that it diverges to
$+\infty$ when $|h_x-h_y|\to\infty$: the potential has the effect of
``flattening'' the interface and suppressing wild fluctuations, in
agreement with the observed macroscopic flatness of phase
boundaries. (Much more stringent conditions have to be imposed on $V$
to actually prove any result.) Note also that the measure
$\pi_{\Lambda_\epsilon}$ depends on the boundary height
$\bar h_\cdot$: if $\bar h_\cdot$ is  fixed so that the
average slope is $\rho\in \mathbb R^d$,
i.e. $\pi_{\Lambda_\epsilon}(h_x-h_{x+e_i})=\rho_i,i\le d$, then we
write $\pi_{\Lambda_\epsilon,\rho}$.

There are various choices of Markov dynamics that admit \eqref{eq:BG}
as stationary reversible measure. A popular choice is the heat-bath or
Glauber dynamics: with rate $1$, independently, each height $h_x(t)$
is refreshed and the new value is chosen from the stationary measure
$\pi_{\Lambda_\epsilon}$ \emph{conditioned} on the values of $h_y(t)$
with $y$ ranging over the nearest neighbors of $x$. Another natural
choice, when the heights are in $\mathbb R$ rather than in
$\mathbb Z$, is a Langevin-type dynamics where each $h_x(t)$ is
subject to an independent Brownian noise, plus a drift that depends on
the height differences between $x$ and its neighboring sites, chosen
so that \eqref{eq:BG} is reversible.

As we mentioned, under reasonable assumptions, a diffusive hydrodynamic limit is expected:
\begin{eqnarray}
  \label{eq:h2}
 \epsilon h_{\epsilon^{-1} x}(\epsilon^{-2}t)\stackrel[\epsilon\to0]{\mathbb P}{\longrightarrow} \phi(x,t)
\end{eqnarray}
where $\phi$ is deterministic. Due to the diffusive scaling of time, the PDE solved by $\phi$ will be of second order and in general non-linear:
\begin{eqnarray}
  \label{eq:h20}
  \partial_t \phi(x,t)=\mu(\nabla \phi(x,t))\sum_{i,j=1}^d \sigma_{i,j}(\nabla \phi(x,t))\frac{\partial^2}{\partial_{x_i}\partial_{ x_j}}\phi(x,t).
\end{eqnarray}
The factors $\mu$ and $\sigma_{i,j}$ have a very different origin,
which is why we have not written the equation in terms of the
combination $\tilde \sigma_{i,j}:=\mu \sigma_{i,j}$ instead.  The
slope-dependent prefactor $\mu>0$ is called \emph{mobility} and will
be discussed in a moment. As for $\sigma_{i,j}$, let the convex function
$\sigma:\rho\in\mathbb R^d\mapsto \sigma(\rho)\in \mathbb R$ denote
the surface tension of the model at slope $\rho$ \cite{Funaki}, i.e. minus the limit as $\epsilon\to0$ of $1/{|\Lambda_\epsilon|}$ times the logarithm of the normalization constant of the probability measure \eqref{eq:BG} when $\pi_{\Lambda_\epsilon}=\pi_{\Lambda_{\epsilon,\rho}}$. Then, $\sigma_{i,j}$
denotes the second derivative of $\sigma$ w.r.t. the $i^{th}$ and
$j^{th}$ argument. Convexity of $\sigma$ implies that the matrix
$\{\sigma_{i,j}(\nabla\phi)\}_{i,j=1,\dots,d}$ is positive definite, so the PDE \eqref{eq:h20} is of parabolic type.
We emphasize that the surface tension, hence $\sigma_{i,j}$, are defined purely in terms of the stationary measure \eqref{eq:BG}. All the dependence on  the Markov dynamics is in the mobility $\mu$.
Remark also that one can rewrite \eqref{eq:h20} in the following more evocative form:
\begin{eqnarray}
  \label{eq:h21}
   \partial_t \phi(x,t)=-\mu(\nabla\phi(x,t))\frac{\delta F[\phi(\cdot,t)]}{\delta \phi(x,t)}
\end{eqnarray}
where $F[\phi(\cdot)]=\int dx\, \sigma(\nabla \phi)$ is the surface
tension functional and $\delta F/\delta\phi$ denotes its first
variation. In other words, the hydrodynamic equation is nothing but
the gradient flow w.r.t. the surface tension functional, modulated by
a slope-dependent mobility prefactor.

Via linear response theory one can guess a Green-Kubo-type expression
for the mobility \cite{SpohnJSP}. This turns out to be given
as\footnote{One can express $\mu$ also via a variational principle,
  see \cite{Spohnbook}.} (say that the heights $h_x$ are discrete, so
that the dynamics is a Markov jump process; the formula is analogous
for Langevin-type dynamics)
\begin{eqnarray}
  \label{eq:GK}
  \mu(\rho)&=&\lim_{\epsilon\to0}\frac1{2|\Lambda_\epsilon|}\pi_{\Lambda_\epsilon,\rho}\left[
  \sum_{x\in\Lambda_\epsilon}\sum_n c_x^n(h)n^2\right]\\
\label{eq:GKriga2} & -&\int_0^\infty dt \lim_{\epsilon\to0}\frac1{|\Lambda_\epsilon|}\sum_{x,x'\in\Lambda_\epsilon}\sum_{n,n'}\mathbb E_{\pi_{\Lambda_\epsilon,\rho}}\left[c_x^n(h(0))n\; c_{x'}^{n'}(h(t))n'\right]
\end{eqnarray}
where $c_x^n(h)$ is the rate at which the height at $x$ increases by
$n\in\mathbb Z$ in configuration $h$, $\mathbb E_{\pi_{\Lambda_\epsilon,\rho}}$
denotes expectation w.r.t. the stationary process started from
$\pi_{\Lambda_\epsilon,\rho}$ and $h(t)$ denotes the configuration at
time $t$.  Note that the first term involves only equilibrium
correlation functions in the infinite volume stationary measure
$\pi_\rho=\lim_{\Lambda\to\mathbb
  Z^d}\pi_{\Lambda_\epsilon,\rho}$\footnote{For models in dimension
  $d\le 2$ the law of the interface does not have a limit as
  $\epsilon\to0$, since the variance of $h_x$ diverges as
  $\epsilon\to0$. However, the law of the \emph{gradients} of $h$ does
  have a limit and the transition rates $c_x^n(h)$ are actually
  functions of the gradients of $h$ only, by translation invariance in the vertical direction. }. The same is not true for
the second one, which involves a time integral of correlations at
different times for the stationary process. These are usually not
explicitly computable even when $\pi_\rho$ is known. It may however
happen for certain models that, by a discrete summation by parts
w.r.t. the $x$ variable, $\sum_x\sum_n c_x^n(h)$ is deterministically
zero, for any configuration $h$: one says then that a \emph{gradient
  condition} is satisfied (a classical example is 
symmetric simple exclusion). In this case \eqref{eq:GKriga2}
identically vanishes and one is in a much better position to prove
convergence to the hydrodynamic equation.

For the ``Ginzburg-Landau (GL)'' model \cite{SpohnJSP} where heights
$h_x$ are continuous variables and the dynamics is of Langevin type, if the potential $V(\cdot)$ is convex and
symmetric then, in any dimension $d$,
the gradient condition is satisfied and moreover the remaining average
in the Green-Kubo formula is immediately computed, leading to a constant
mobility: $\mu(\rho)=1$.  In this situation, Funaki and Spohn \cite{FuSpo}  proved
convergence of the height profile to (the weak solution of) \eqref{eq:h20} for the GL
model, for any $d\ge1$. (They look at weak solutions because for $d>1$ it is not known whether
the surface tension of the GL model  is $C^2$ and the coefficients
$\sigma_{i,j}$ are well defined and smooth). Until recently, to my knowledge,
there was no other known interface model in dimension $d>1$ where 
 mathematical results of this type were available.

Before presenting our recent results for $(2+1)$-dimensional interface
dynamics let us make two important observations:
\begin{itemize}
\item Not only for most natural interface dynamics in dimension $d>1$
  one is unable to prove a hydrodynamic convergence of the type
  \eqref{eq:h2}: the situation is actually much worse. As
  \eqref{eq:h2} suggests, the correct time-scale for the system to
  reach stationarity (measured either by
  $T_{rel}:=1/{\rm gap}(\mathcal L)$, with gap$(\mathcal L)$ denoting
  the spectral gap of the generator, or by the so-called total
  variation mixing time $T_{mix}$) should be of order $\epsilon^{-2}$
  (logarithmic corrections are to be expected for the mixing time). On
  the other hand, for most natural models it not even proven that such
  characteristic times are upper bounded by a polynomial of
  $\epsilon^{-1}$! For instance, for the well-known
  $(2+1)$-dimensional SOS model at low temperature, the best known
  upper bound for $T_{rel}$ and $T_{mix}$ is a rather poor
  $O(\exp(\epsilon^{-1/2+o(1)}))$ \cite[Th. 3]{CLMST}.
  
\item In dimension $d=1$, natural Markov dynamics of discrete
  interfaces are provided by conservative lattice gases on $\mathbb Z$
  (e.g. symmetric exclusion processes or zero-range processes), just
  by interpreting the number of particles at site $x$ as the interface
  gradient $h_x-h_{x-1}$ at $x$. Similarly, conservative
  continuous spin models on $\mathbb Z$ translate into Markov dynamics
  for one-dimensional interface models with continuous heights. Then,
  a hydrodynamic limit for the height function follows from that for
  the particle density (see e.g. \cite[Ch. 4 and 5]{KL} for the
  symmetric simple exclusion and for a class of zero-range processes,
  and for instance \cite{Fritz} for the $d=1$ Ginzburg-Landau
  model). For $d>1$, instead, there is in general no obvious way of
  associating a height function to a particle system on $\mathbb
  Z^d$. Also, for $d=1$ there are robust methods to prove that the inverse
  spectral gap is
  $T_{rel}=O(\epsilon^{-2})$, see e.g. \cite{KL,Caputo}.
\end{itemize}

\subsection{Reversible tiling dynamics, mixing time and  hydrodynamic equation}

\begin{figure}
\begin{center}
\includegraphics[height=2cm]{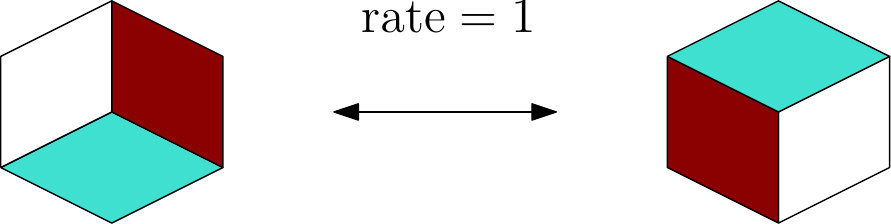}
\caption{The elementary updates for the Glauber dynamics of the
  lozenge tiling model correspond to the rotation of three lozenges
  (equivalently, three dimers) around a hexagonal face. The transition
  rate is $1$ both for the update and the reverse one. Note that in
  the three-dimensional corner growth model only the forward
  transition would be allowed.}
\label{fig:flips}
\end{center}
\end{figure}

In this section, I briefly review a series of results obtained in
recent years in collaboration with Pietro Caputo, Beno\^it Laslier and
Fabio Martinelli.  In these works we study $(2+1)$-dimensional
interface dynamics where the height function
$\{h_x\}_{x\in\Lambda_\epsilon}$ is discrete and is given by the
height function of a tiling model, either by lozenges or by dominoes,
as explained in Section \ref{sec:newAKPZ}. In contrast with the
$(2+1)$-dimensional Anisotropic KPZ growth models described in Section
\ref{sec:newAKPZ}, that are also Markov dynamics of tiling models,
here we want a \emph{reversible} process because we wish to model interface
evolution at thermal equilibrium.  A natural candidate is the
``Glauber'' dynamics obtained by giving rate $1$ to the elementary
rotations of tiles around faces of the graph, see
Fig. \ref{fig:flips} for the case of lozenge tilings.

In terms of the height function, elementary moves correspond to
changing the height by $\pm1$ at single sites. Since all elementary
rotations have the same rate, the uniform measure over the finitely
many tiling configurations in $\Lambda_\epsilon$ is reversible. As a
side remark, this measure can be written in the Boltzmann-Gibbs form
\eqref{eq:BG} with a potential $V$ taking values $0$ or $+\infty$.
Let us also remark that, as discussed in \cite{CMT}, this dynamics is
equivalent to the zero-temperature Glauber dynamics of the
\emph{three-dimensional} Ising model with Dobrushin boundary
conditions.

In agreement with the discussion of the previous section, if the tiled
region is a reasonably-shaped domain $\Lambda_\epsilon$ of diameter
$O(\epsilon^{-1})$, one expects $T_{rel}$ and $T_{mix}$ to be $\approx \epsilon^{-2}$ and the height profile to converge
under diffusive rescaling to the solution of a parabolic PDE. Until
 recently, however, all what was known rigorously was that
$T_{rel}$ and $T_{mix}$ are upper bounded as
$O(\epsilon^{-n})$ for some finite $n>2$!
\begin{op}
  This polynomial upper bound was proven in \cite{LRS} for the Glauber
  dynamics on either lozenge or domino tiling (the same
  proof works for tilings associated to the dimer model on  certain graphs with both hexagonal and square faces,
  as shown in \cite{LasTonDomino}). The method does not seem to work,
  however, for general planar bipartite graphs. For instance, a
  polynomial upper bound for $T_{rel}$ or $T_{mix}$ of the Glauber
  dynamics of the dimer model on the square-octagon graph (see Fig. 9 in \cite{Kenyon})
  is still unproven. \end{op}

Under suitable conditions,  we improved this
$O(\epsilon^{-n})$ upper bound into an almost optimal one:
\begin{thm}[Informal statement]
If the boundary height on
$\partial\Lambda_\epsilon$ is such that the average height under the measure
$\pi_{\Lambda_\epsilon}$ tends to an affine function as
$\epsilon\to0$, then $T_{mix}$ and $T_{rel}$ are $O(\epsilon^{-2+o(1)})$
\cite{CMT,LasTonDomino}.
\end{thm}
 Later \cite{LasTonMacro}, we proved a result in the
same spirit under the sole assumption that the limit average height
profile is smooth and in particular has no ``frozen regions'' \cite{Kenyon}.

Let us emphasize that there are very natural domains
$\Lambda_\epsilon$ such that the average equilibrium height profile in
the $\epsilon\to0$ limit does have ``frozen regions'':
 \begin{op}
   Let $\Lambda_\epsilon$ be the hexagonal domain of side
   $\epsilon^{-1}$, see Fig.  \ref{fig:hexa}.  Prove that, for the
   lozenge tiling Glauber dynamics, $T_{rel}$ and/or $T_{mix}$ are
   $O(\epsilon^{-2+o(1)})$. The best upper bound that can be extracted
   from \cite{Wilson} plus the so-called Peres-Winkler censoring
   inequalities \cite{PW} is $O(\epsilon^{-4+o(1)})$.
 \end{op}

The proofs of the previously known polynomial upper bounds on the
mixing time were based on smart and rather simple path coupling
arguments \cite{LRS}.  To get our almost optimal bounds  \cite{CMT,LasTonDomino,LasTonMacro}, there are at
least two new inputs:
 \begin{itemize}
 \item our proof consists in a comparison between the actual
   interface dynamics and an auxiliary one that evolves on
   almost-diffusive time-scales $\approx \epsilon^{-2+o(1)}$ and that
   essentially follows the conjectural hydrodynamic motion where
   interface drift is proportional to   its curvature;

  
 \item to control the auxiliary process, we crucially need very  refined estimates on  height fluctuation for the equilibrium measure $\pi_{\Lambda_\ell}$ on domains of mesoscopic size $\ell=\epsilon^{-a},1/2\le a\le1$, with various types of boundary conditions.
 \end{itemize}

\begin{figure}
\begin{center}
\includegraphics[height=3cm]{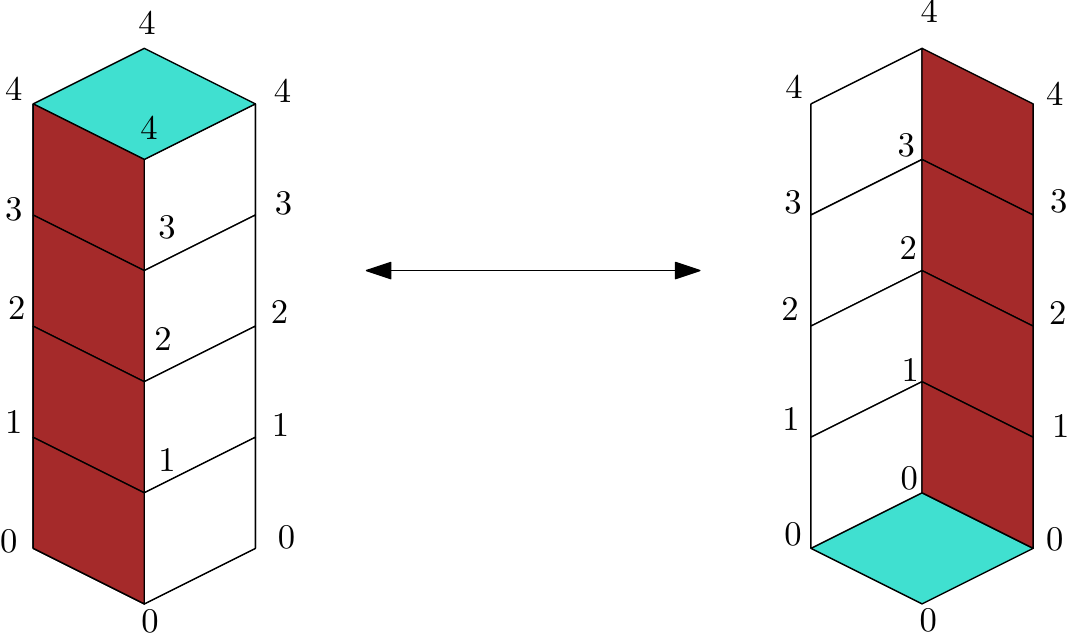}
\caption{A ``tower move'' transition  with $n=4$ and the reverse transition. The transition rates equal $1/n=1/4$. Note that the height decreases by $1$ at four points.}
\label{fig:flippone}
\end{center}
\end{figure}

For the Glauber dynamics with elementary moves as in Fig. \ref{fig:flips}, it
seems hopeless to prove a hydrodynamic limit on the diffusive
scale. In particular, no form of ``gradient condition'' is
satisfied. Fortunately, there exists a more friendly variant of
the Glauber dynamics, introduced in \cite{LRS}, where a single update
consists in ``tower moves'' changing the height by the same amount $\pm1$ at $n\ge0$
aligned sites, as in Fig. \ref{fig:flippone}.
The integer $n$ is not fixed here, in
fact transitions with any $n$ are allowed but the transition rate
decreases with $n$ and actually it is taken to equal to $1/n$. It is
immediate to verify that this dynamics is still reversible w.r.t. the
uniform measure.  For this modified dynamics, together with B. Laslier
we realized in \cite{LasTonHydro0} that a microscopic summation by
parts implies that the term \eqref{eq:GKriga2} in the definition of
the mobility vanishes, and actually we could explicitly compute $\mu$, that turns out to be
\emph{non-trivial and non-linear}:
\begin{eqnarray}
  \label{eq:munuova}
  \mu(\rho)=\frac1\pi\frac{\sin(\pi \rho_1)\sin(\pi\rho_2)}{\sin(\pi(\rho_1+\rho_2))}.
\end{eqnarray}
Recall that, in contrast, the mobility of the Ginzburg-Landau model is
slope-independent \cite{SpohnJSP}.  Later, in \cite[Th. 2.7]{LasTonHydro}, we
could turn our arguments into a full proof of convergence of the
height profile to the solution of the PDE:
\begin{thm}[Informal statement] 
  If the initial profile $\phi_0$ is sufficiently smooth, one has for every $t>0$
  \begin{eqnarray}
    \lim_{\epsilon\to0} \frac1{|\Lambda_\epsilon|}\sum_{x\in\Lambda_\epsilon}\mathbb E\left|\epsilon h_x(\epsilon^{-2}t)-\phi(\epsilon x,t)\right|^2=0,
  \end{eqnarray}
with $\phi(x,t)$ the solution of \eqref{eq:h20}
\end{thm}
 (For technical reasons, we
had to work with periodic instead of Dirichlet boundary conditions).
A couple of comments are in order:
\begin{itemize}
\item As the reader may have noticed, the function \eqref{eq:munuova}
  is exactly the same as the ``speed function'' $v(\rho)$ of the
  growth model discussed in Section \ref{sec:newAKPZ}, see formula
  \eqref{eq:vmu}.  This is not a mere coincidence. Actually, one may see
  this equality as an instance of the so-called Einstein relation
  between diffusion and conductivity coefficients \cite{Spohnbook}.
  
\item We mentioned earlier that convergence of the height profile of
  the Ginzburg-Landau model to the limit PDE has been proved
  \cite{FuSpo} only in a weak sense. In our case, instead, we have
  strong convergence to classical solutions of \eqref{eq:h20} that
  exist globally because the coefficients
  $\mu(\cdot),\sigma_{i,j}(\cdot)$ turn out to be smooth functions of
  the slope.\footnote{ The apparent singularity of the formula
    \eqref{eq:munuova} for $\mu(\cdot)$ when $\rho_1+\rho_2=0$ is not
    really dangerous: recall from Section \ref{sec:newAKPZ} that the slope $\rho$ is constrained in
    the triangle
    $\mathbb T=\{(\rho_1,\rho_2):0\le \rho_1,\rho_2, \rho_1+\rho_2\le
    1\}$ so that the mobility is $C^\infty$ and strictly positive in the interior of
    $\mathbb T$.}
  
\item A fact that plays a crucial role in the proof of the
  hydrodynamic limit is that the PDE \eqref{eq:h20} contracts the
  $L^2$ distance $D_2(t)=\int dx(\phi^{(1)}(x,t)-\phi^{(2)}(x,t))^2$
  between solutions. I believe this is not a trivial or general fact:
  in fact, to prove contraction \cite{LasTonHydro}, we use the
  specific form \eqref{eq:munuova} of $\mu$ and the explicit
  expression of $\sigma_{i,j}$ for the dimer model. (Note that if the
  mobility were constant, as it is for the Ginzburg-Landau model,
  $L^2$ contraction would just be a consequence of convexity of the
  surface tension $\sigma$).  I think it is an intriguing question to
  understand whether the identities (see
  \cite[Eqs. (6.19)-(6.22)]{LasTonHydro}) leading to
  $d D_2(t)/dt\le 0$ have any thermodynamic interpretation.
 \end{itemize}

 To conclude this review, let us mention that new dynamical phenomena,
 taking place on time-scales much longer than diffusive, can occur at
 low temperature, for interface models undergoing a so-called
 ``roughening transition''. That is, up to now we considered
 situations where the equilibrium Gibbs measure for the interface in a
 $\epsilon^{-1}\times \epsilon^{-1}$ box $\Lambda_\epsilon$ scales to
 a massless Gaussian field as $\epsilon\to0$ and in particular
 ${\rm Var}_{\pi_{\Lambda_\epsilon}}(h_x)\approx \log(1/\epsilon)$ if
 $x$ is, say, the center of the box. The interface is said to be
 ``rough'' in this case, because fluctuations diverge as $\epsilon\to0$. For some interface models, notably the well-known
 Solid-on-Solid (SOS) model where the potential $V$ in \eqref{eq:BG}
 equals $T^{-1} |h_x-h_y|$ and heights are integer-valued and fixed to $0$ around the
 boundary, it is known that at low enough temperature $T$ the
 interface is instead rigid, with
 $\limsup_{\epsilon\to0} {\rm
   Var}_{\pi_{\Lambda_\epsilon}}(h_x)<\infty$, while the variance
 grows logarithmically at high temperature \cite{FS}. The temperature
 $T_r$ separating these two regimes is called ``roughening
 temperature''.

In a work with P. Caputo, E. Lubetzky, F. Martinelli and A. Sly \cite{CLMST} we discovered that rigidity of the interface can 
produce a dramatic slowdown of the dynamics, if the interface is constrained to stay above a fixed level, say level $0$:
\begin{thm}\cite{CLMST} Consider the Glauber dynamics for the $(2+1)$-dimensional SOS model at low enough temperature, 
with $0$ boundary conditions on $\partial \Lambda_\epsilon$ and with the positivity constraint  $h_x\ge0 $ for every $x\in\Lambda_\epsilon$. Then, the relaxation and mixing times satisfy
\begin{eqnarray}
  T_{mix}\ge T_{rel}\ge c \exp[{c\;\epsilon^{-1}}]
\end{eqnarray}
for some positive, temperature-dependent constant $c$.
\end{thm}
What we actually prove is that there is a cascade of metastable
transitions, occurring on all time-scales $\exp({\epsilon^{-a}})$,
$a<0\le 1$.  Strange as this may look, these results \emph{do not}
exclude that a hydrodynamic limit on the \emph{diffusive scale}, as in
\eqref{eq:h2}-\eqref{eq:h20},  might occur. That is, the rescaled
height profile $\epsilon h_{\epsilon^{-1}x}(\epsilon^{-2}t)$ could
follow an equation like \eqref{eq:h20}, so that at times $\gg \epsilon^{-2}$ the
profile would be macroscopically zero (because $\phi\equiv 0$ is the
equilibrium point of the PDE \eqref{eq:h20} with zero boundary conditions) but smaller-scale height
fluctuations would need enormously more time, of the order $T_{mix}\approx \exp(\epsilon^{-1})$, to
relax to equilibrium.

We are light years away from being able to actually prove a hydrodynamic limit for the $(2+1)$-dimensional SOS model.
The following open problem is given just to show how little we know in this respect:
\begin{op}
  Take the Glauber dynamics for the $(2+1)$-d SOS model at low temperature, with
  initial condition $\epsilon\, h_x=1$ for every
  $x\in\Lambda_\epsilon$. Is it true that, for some $N<\infty$, at
  time $t=1/\epsilon^N$ all rescaled heights $\epsilon \,h_x(t)$ are
  with high probability lower than, say, $1/2$ (which is much larger
  than $\epsilon\log\epsilon^{-1}$, that is the typical value under the
  equilibrium measure of
  $\max_{x\in\Lambda_\epsilon}[\epsilon \,h_x]$)?
\end{op}

\subsubsection*{Acknowledgements}
This work  was partially funded by the ANR-15-CE40-0020-03 Grant LSD, by the CNRS
PICS grant ``Interfaces al\'eatoires discr\`etes et dynamiques de
Glauber'' and by MIT-France Seed Fund ``Two-dimensional Interface
Growth and Anisotropic KPZ Equation''.

\end{document}